\documentclass[generic,11pt,preprint]{imsart}
\usepackage{lineno}
\modulolinenumbers[5]

\usepackage{url}
\usepackage{verbatim}
\RequirePackage[OT1]{fontenc}
\RequirePackage[numbers]{natbib}
\RequirePackage{hypernat}

\usepackage{mathscinet}
\usepackage[utf8]{inputenc} 
\usepackage{textcomp} 
\usepackage{graphicx}  
\usepackage{flafter}  
\usepackage{natbib} 
\usepackage{amsmath,amssymb,amsthm}  
\usepackage{bm}  
\usepackage[pdftex,bookmarks,colorlinks,breaklinks]{hyperref}  
\usepackage{url}
\usepackage{subfigure}
\usepackage{verbatim}
\usepackage{wrapfig}
\usepackage{enumerate}
\hypersetup{linkcolor=red,citecolor=blue,filecolor=cyan,urlcolor=blue} 
\numberwithin{equation}{section}
\theoremstyle{plain}
\newtheorem{definition}{Definition}[section]
\newtheorem{theorem}{Theorem}[section]
\newtheorem{corollary}{Corollary}[section]

\newtheorem{lemma}{Lemma}[section]

\newtheorem{remark}{Remark}
\newtheorem{fact}{Fact}

\newcommand{\lmax}[1]{\lambda_{\max}\left(#1\right)}
\newcommand{\lmin}[1]{\lambda_{\min}\left(#1\right)}
\newcommand{\prob}[1]{\mathbb{P}\left(#1\right)}
\newcommand{\dotp}[2]{\left\langle#1,#2\right\rangle}
\newcommand{\m}{\mathcal}
\newcommand{\mb}{\mathbb}

\newcommand{\tr}{\mbox{tr}\,}
\def\r{\right}
\def\l{\left}

\newcommand{\mx}{\mbox{\footnotesize{max}\,}}
\newcommand{\mn}{\mbox{\footnotesize{min}\,}}

\begin{document}

\begin{frontmatter}
\title{On Some Extensions of Bernstein's Inequality for Self-Adjoint Operators}
\begin{aug}
\author{\fnms{Stanislav} \snm{Minsker}\thanksref{t2}\ead[label=e1]{minsker@usc.edu}}

 \thankstext{t2}{Department of Mathematics, University of Southern California}
\runauthor{S. Minsker}

\affiliation{University of Southern California}
\printead{e1}\\
\end{aug}

\begin{abstract}
We present some extensions of Bernstein's concentration inequality for random matrices. 
This inequality has become a useful and powerful tool for many problems in statistics, signal processing and theoretical computer science. 
The main feature of our bounds is that, unlike the majority of previous related results, they do not depend on the dimension $d$ of the ambient space. 
Instead, the dimension factor is replaced by the ``effective rank'' associated with the underlying distribution that is bounded from above by $d$. 
In particular, this makes an extension to the infinite-dimensional setting possible. 
Our inequalities refine earlier results in this direction obtained by D. Hsu, S. M. Kakade and T. Zhang.  
\end{abstract}

\begin{keyword}
Bernstein's inequality \sep Random Matrix \sep Effective rank \sep Concentration inequality \sep Large deviations
\end{keyword}

\end{frontmatter}


\section{Introduction}

Theoretical analysis of many problems in applied mathematics can often be reduced to problems about random matrices. 
Examples include numerical linear algebra (randomized matrix decompositions, approximate matrix multiplication), statistics and machine learning (low-rank matrix recovery, covariance estimation), mathematical signal processing, among others. 

Often, resulting questions can be reduced to estimates for the expectation $\mb E\l\| \sum\limits_{i=1}^n X_i - \mb EX_i \r\|$ or probability $\prob{\left\|\sum\limits_{i=1}^n X_i - \mb EX_i\right\|>t}$, where $\left\{X_i \right\}_{i=1}^n$ is a finite sequence of random matrices and $\|\cdot\|$ is the operator norm. 
Some of the early work on the subject includes the papers on noncommutative Khintchine inequalities \cite{lust1991non,lust1986inegalites}; these tools were used by M. Rudelson \cite{rudelson1999random} to establish the bounds for the deviations of a sample covariance matrix from its population counterpart. 
The pioneering work of R. Ahlswede and A. Winter \cite{Ahlswede2002Strong-converse00} demonstrated how to adapt the Laplace transform technique the non-commutative framework. 
Ahlswede and Winter's work relied on Golden-Thompson inequality \cite{sidne1965lower,thompson1965inequality};
a different approach that often yields sharper results was developed in the works by R.I. Oliveira and J. Tropp \cite{Oliveira2010Concentration-o00,tropp2012user,tropp2} and led to the versions of Chernoff, Bernstein and Friedman bounds for random matrices. 

One of the difference in the tail bounds between the scalar and matrix versions of these concentration inequalities is the multiplicative dimension factor appearing in the bounds. 
For instance, the matrix Bernstein inequality \cite{tropp2012user} states that, given a sequence $X_1,\ldots,X_n\in \mb C^{d\times d}$ of i.i.d. self-adjoint copies of a random matrix $X$ such that $\|X\|\leq U$ almost surely, $\prob{\l\| \sum_{i=1}^n (X_i -\mb EX)  \r\|\geq t}\leq 2d \exp\l(-\frac{t^2/2}{\sigma^2+Ut/3} \r)$ for any $t>0$; here, $\sigma^2$ is a ``variance parameter'' defined below. 
While unavoidable in general, factor $d$ can often be replaced by the ``effective rank'' that is always bounded above by $d$ but can be much smaller in some cases. 
Obtaining the bounds that depend on effective rank is the main purpose of this paper. 
Results in this direction were previously presented in \cite{Hsu2012Tail-inequaliti00}: for instance, in the i.i.d. case authors have shown that 
$\prob{\l\| \sum_{i=1}^n (X_i -\mb EX)  \r\|\geq \sigma\sqrt{\frac{2t}{n}} + \frac{Ut}{3n}}\leq 2 r \frac{t}{e^t - t - 1}$, where $r$ is the effective rank defined below. 
We show that in many cases the tail behavior can be further improved through the modification of original Tropp's method \cite{tropp2012user} based on Lieb's concavity theorem: for instance, factor $t$ appearing in the numerator on the right-hand side of the previous inequality can be removed. 

The paper is organized as follows: section \ref{sec:prelim} introduces main notation and background material, section \ref{sec:main} presents the main results, and section \ref{sec:applications} contains examples. 


\section{Preliminaries}
\label{sec:prelim}

In this section, we introduce notation and recall several useful facts from linear algebra, matrix analysis and probability theory that we rely on in the subsequent exposition. 

\subsection{Definitions and notation}

Given $A\in \mb C^{d_1\times d_2}$, let $A^\ast\in \mb C^{d_2\times d_1}$ be the Hermitian adjoint of $A$. 
If $A$ is self-adjoint, we will write $\lambda_{\mx}(A)$ and $\lambda_{\mn}(A)$ for the largest and smallest eigenvalues of $A$. 
Next, we will introduce the matrix norms used in the paper. 
Everywhere below, $\|\cdot\|$ stands for the operator norm $\|A\|:=\sqrt{\lambda_{\mx}(A^\ast A)}$. 
If $d_1=d_2=d$, we denote by $\tr A$ the trace of $A$.
The nuclear norm $\|A\|_1$ of $A\in \mb C^{d_1\times d_2}$ is defined as 
$\|A\|_1=\tr(\sqrt{A^*A})$, where $\sqrt{A^*A}$ is a nonnegative definite matrix such that $(\sqrt{A^*A})^2=A^\ast A$. 
The Frobenius (or Hilbert-Schmidt) norm is $\|A\|_{\mathrm{F}}=\sqrt{\tr(A^\ast A)}$, and the associated inner product is 
$\dotp{A_1}{A_2}=\tr(A_1^\ast A_2)$. 
For $Y\in \mb C^d$, $\l\| Y \r\|_2$ stands for the standard Euclidean norm of $Y$. 

Given two self-adjoint operators $A$ and $B$, we will write $A\succeq B \ (\text{or }A\succ B)$ iff $A-B$ is nonnegative (or positive) definite. 

For a sequence of random matrices $X_1,\ldots,X_n$, $\mb E_j[\,\cdot \,]$ will stand for the conditional expectation 
$\mb E[\,\cdot\,| X_1,\ldots,X_{j}]$, where expectation of a matrix is evaluated element-wise.
Finally, for $a,b\in \mb R$, set $a\vee b:=\max(a,b)$ and $a\wedge b:=\min(a,b)$. 

Next, we introduce the notion of ``effective rank'' \cite{koltchinskii2016asymptotics}. 
Let $A\in \mb C^{d\times d}$ be nonnegative definite. 
The effective rank of $A$ is defined as 
\begin{align}
\label{eq:effdim}
r(A):=\frac{\tr \l(A \r) }{\l\| A  \r\|}.
\end{align}
Clearly, $r(A)\leq d$, and $r(A)$ can be much smaller than $d$ if $A$ has many eigenvalues that are close to $0$ in absolute value. 
Finally, we recall the definition of a function of a matrix-valued argument that will be frequently used in the paper:
\begin{definition}
Given a real-valued function $f$ defined on an interval $\mb T\subseteq \mb R$ and a self-adjoint $A\in \mb C^{d\times d}$ with the eigenvalue decomposition 
$A=U\Lambda U^\ast$ such that $\lambda_j(A)\in \mb T,\ j=1,\ldots,d$, define $f(A)$ as 
$f(A)=Uf(\Lambda) U^\ast$, where 
\[
f(\Lambda)=f\l( \begin{pmatrix}
\lambda_1 & \,  & \,\\
\, & \ddots & \, \\
\, & \, & \lambda_d
\end{pmatrix} \r)
=\begin{pmatrix}
f(\lambda_1) & \,  & \,\\
\, & \ddots & \, \\
\, & \, & f(\lambda_d)
\end{pmatrix}.
\] 
\end{definition}

\subsection{Tools from linear algebra}

We recall several facts from linear algebra, matrix analysis and probability theory that play important roles in our arguments. 
\begin{fact}
\label{fact:01}
Let $A\in \mb C^{d\times d}$ be a self-adjoint matrix, and $f_1, \ f_2$ be two real-valued functions such that $f_1(\lambda_j)\geq f_2(\lambda_j)$ for $j=1,\ldots,d$. 
Then $f_1(A)\succeq f_2(A)$. 
\end{fact}

\begin{fact}
\label{fact:02}
Let $A,B\in \mb C^{d\times d}$ be two self-adjoint matrices such that $A\succeq B$. 
Then $\lambda_j(A)\geq \lambda_j(B), \ j=1,\ldots, d$, where $\lambda_j(\cdot)$ stands for the $j$-th largest eigenvalue.   
Moreover, $\tr (e^A)\geq \tr (e^B)$.
\end{fact}

\begin{fact}
\label{fact:03}
Matrix logarithm is operator monotone: if $A\succ 0, \ B\succ 0$ and $A\succeq B$, then $\log(A)\succeq \log(B)$. 
\end{fact}
\begin{proof}
See \cite{bhatia1997matrix}. 
\end{proof}

\begin{fact}
\label{fact:04}
Let $A, \ B\in \mb C^{d\times d}$ be a self-adjoint matrices such that $A\succeq B$. 
Then for any $Q\in \mb C^{d\times d}$, $Q A Q^\ast \succeq Q B Q^\ast$.
\end{fact}

\begin{fact}[Lieb's concavity theorem]
\label{fact:05}
Given a fixed self-adjoint matrix $H$, the function
\[
A\mapsto \tr\exp\l(H+\log(A)\r)
\]
is concave on the cone of positive definite matrices.
\end{fact}
\begin{proof}
See \cite{lieb1}. 
\end{proof}

\begin{fact}
\label{fact:06}
Let $f:\mb R\mapsto \mb R$ be a convex function. 
Then $A\mapsto \tr f(A)$ is convex on the set of self-adjoint matrices. 
In particular, for any self-adjoint matrices $A,B$, 
\[
\tr f\l( \frac{A+B}{2} \r)\leq \frac{1}{2}\tr f(A) + \frac{1}{2}\tr f(B).
\] 
\end{fact}
\begin{proof}
This is a consequence of Peierls inequality, see Theorem 2.9 in \cite{carlen2010trace} and the comments following it. 
\end{proof}

Finally, we introduce the Hermitian dilation which allows to reduce many problems involving general rectangular matrices to the case of Hermitian operators. 
Given the rectangular matrix $A\in\mb C^{d_1\times d_2}$, the Hermitian dilation $\m H: \mb C^{d_1\times d_2}\mapsto \mb C^{(d_1+d_2)\times (d_1+d_2)}$ is defined as
\begin{align}
\label{eq:dilation}
&
\m H(A)=\begin{pmatrix}
0 & A \\
A^\ast & 0
\end{pmatrix}.
\end{align}
Since 
$\m H(A)^2=\begin{pmatrix}
A A^\ast & 0 \\
0 & A^\ast A
\end{pmatrix},$ 
it is easy to see that $\| \m H(A) \|=\|A\|$.

\section{Main results}
\label{sec:main}

We are ready to state and prove our main results. 
We begin by stating a version of Bernstein's inequality for self-adjoint random matrices, and will later deduce the result for rectangular matrices from it.

\begin{theorem}
\label{th:independent}
Let $X_1,\ldots, X_n\in \mb C^{d\times d}$ be a sequence of independent self-adjoint random matrices such that $\mb EX_i=0, \ i=1,\ldots,n$ and 
$\sigma^2\geq\l\|\sum\limits_{i=1}^n \mb E X_i^2\r\|$. 
Assume that $\|X_i\|\leq U$ almost surely for all $1\leq i\leq n$ and some positive $U\in \mb R$. 
Then, for any $t\geq \frac{1}{6}\l( U+\sqrt{U^2 + 36\sigma^2} \r)$, 
\begin{align}
\label{eq:a05}
&
\prob{\left\|\sum_{i=1}^n X_i\right\|>t}\leq 14 \, r\left(\sum\limits_{i=1}^n \mb EX_i^2\right)
\exp\left[- \frac{t^2/2}{\sigma^2+tU/3} \right],
\end{align}
where $r(\cdot)$ stands for the effective rank defined in \eqref{eq:effdim}.
\end{theorem}
\begin{remark}
\label{remark:1}
Condition $t\geq \frac{1}{6}\l( U+\sqrt{U^2 + 36\sigma^2} \r)$ is not very restrictive as it is equivalent to $t^2\geq \sigma^2 + Ut/3$. 
While the usual Bernstein's inequality holds even when $t^2 < \sigma^2 + Ut/3$, resulting bound is not very useful. 
\end{remark}

\begin{proof}[Proof of Theorem \ref{th:independent}]
Our proof follows the general approach developed in \cite{tropp2012user} to establish the submultiplicativity properties of the trace moment generating function. 
Let $\phi(\theta)=e^{\theta}-\theta-1$. 
Note that $\phi$ is nonnegative and increasing on $[0,\infty)$. 

Let $S_n:=\sum\limits_{i=1}^n X_i$.
Choose any $\theta>0$; the following inequality holds for the largest eigenvalue of $S_n$:
\begin{align}
\label{z2}
\nonumber
\prob{\lmax{S_n}>t}&=\prob{\lmax{\theta S_n}>\theta t} = \prob{\phi \l( \lmax{\theta S_n} \r) > \phi\l( \theta t \r)   } \\ 
&
=\prob{\lmax{\phi(\theta S_n)}>\phi(\theta t)}
\leq \prob{\tr \phi(\theta S_n)>\phi(\theta t)}
\leq \frac{\mb E\tr \, \phi(\theta S_n) }{\phi(\theta t)},
\end{align}
where we have used the fact that $\phi \l( \lmax{\theta S_n} \r) = \lmax{\phi(\theta S_n)}$, as well as Markov's inequality. 
We will next establish the following semidefinite relations:
\begin{lemma}
\label{lemma:semidef}
Under the assumptions of Theorem \ref{th:independent},
\begin{equation}
\label{eq:z1}
 \log \mb E e^{\theta X_i}\preceq \frac{\phi(\theta U)}{U^2}\mb EX_i^2, \ i=1,\ldots,n.
\end{equation}
\end{lemma}
\begin{proof}
See \ref{sec:lemma_pf}.
\end{proof}

\noindent 
Observe that $\exp(\theta S_n) = \exp \l( \theta S_{n-1} + \log e^{\theta X_n} \r)$. 
Since $\mb E S_n = 0$, 
\begin{align*}
\mb E\,\tr \, \phi(\theta S_n)&=\mb E\, \tr \l(\exp(\theta S_{n-1}+\log e^{\theta X_n})-I_d \r)=\\
&
=\mb E\,\mb E_{n-1} \tr \l( \exp \l(\theta S_{n-1}+\log e^{\theta X_n}\r)-I_d \r), 
\end{align*}
where $I_d\in \mb C^{d\times d}$ is the identity matrix. 
Lieb's concavity theorem (fact \ref{fact:05}) with $H:=\theta S_{n-1}$ and $A:=e^{\theta X_n}$, together with Jensen's inequality for conditional expectations imply that
\[
\mb E_{n-1} \tr \l(\exp \l( \theta S_{n-1}+\log e^{\theta X_n}\r)-I_d \r) \leq
\tr \l( \exp \l( \theta S_{n-1}+\log \mb E e^{\theta X_n}\r)-I_d \r),
\]
hence we obtain that $\mb E\,\tr \, \phi(\theta S_n) \leq \mb E \tr \l( \exp \l( \theta S_{n-1}+\log \mb E e^{\theta X_n}\r)-I_d \r)$.
Iterating the argument, we get 
\begin{align*}
\mb E\tr \, \phi(\theta S_n)& 
\leq \mb E \tr \l( \exp \l( \theta S_{n-2}+\log \mb E e^{\theta X_n}+\log \mb E e^{\theta X_{n-1}}\r)-I_d \r) \\
&
\leq \ldots \leq \tr \l( \exp\l(\sum_{i=1}^n \log \mb E e^{\theta X_i}\r)-I_d \r).
\end{align*}
Together with (\ref{eq:z1}), this inequality implies that
\begin{align}
\label{eq:n1}
\mb E\,\tr \phi(\theta S_n)\leq
\tr\l(\exp\l( \frac{\phi(\theta U)}{U^2}\sum_{i=1}^n \mb EX_i^2 \r) - I_d\r).
\end{align}
Note that, since $\mb E S_n^2=\sum_{i=1}^n \mb EX_i^2$ is nonnegative definite, it admits a nonnegative definite square root 
$\l(\mb E S_n^2 \r)^{1/2}$, hence we can write 
\begin{align}
\label{eq:n2}
\nonumber
\exp\Big(&\frac{\phi(\theta U)}{U^2}\mb E S_n^2 \Big)-I_d = \frac{\phi(\theta U)}{U^2}\mb E S_n^2 + \frac{\l( \frac{\phi(\theta U)}{U^2}\mb E S_n^2 \r)^2}{2!}+\ldots+
\frac{\l(\frac{\phi(\theta U)}{U^2}\mb E S_n^2 \r)^k}{k!}+\ldots
\\
\nonumber
&=\frac{\phi(\theta U)}{U^2} \l(\mb E S_n^2\r)^{1/2}\left(1+\frac{1}{2!}\frac{\phi(\theta U)}{U^2}\mb ES_n^2+\ldots
+\frac{1}{n!}\left(\frac{\phi(\theta U)}{U^2}\mb ES_n^2\right)^{n-1}+\ldots\right)\l(\mb E S_n^2\r)^{1/2} \\
\nonumber
&
\preceq \frac{\phi(\theta U)}{U^2} \l(\mb E S_n^2\r)^{1/2}\l(1+\frac{1}{2!}\frac{\phi(\theta U)}{U^2}\|\mb ES_n^2\|+\ldots+
\frac{1}{n!}\left(\frac{\phi(\theta U)}{U^2}\|\mb ES_n^2\|\right)^{n-1}+\ldots\r) \l(\mb E S_n^2\r)^{1/2} \\
&
\preceq \frac{\phi(\theta U)}{U^2} \mb ES_n^2 \, \frac{\exp(\frac{\phi(\theta U)}{U^2}\sigma^2)-1}{\sigma^2 \frac{\phi(\theta U)}{U^2}}
\preceq 
\frac{\mb E S_n^2}{\sigma^2}\exp \l(\frac{\phi(\theta U)}{U^2}\sigma^2 \r),
\end{align}
where we used fact \ref{fact:04} to go from line 2 to line 3. 
Combination of \eqref{z2} with \eqref{eq:n2} and fact \ref{fact:02} (trace monotonicity) yields the inequality
\begin{align}
\label{eq:n10}
& 
\prob{\lmax{S_n}>t}\leq \tr\left(\frac{\mb E S_n^2}{\sigma^2}\right)\frac{\exp\l(\frac{\phi(\theta U)}{U^2}\sigma^2\r)}{\phi(\theta t)} 
= \tr\left(\frac{\mb E S_n^2}{\sigma^2}\right) \exp\l(\frac{\phi(\theta U)}{U^2}\sigma^2 - \theta t\r) \frac{\exp(\theta t)}{\phi(\theta t)}.
\end{align}
Note that, whenever $\theta U <3$,  
\begin{align}
\label{eq:n20}
&
\frac{\phi(\theta U)}{U^2} = \frac{1}{U^2}\sum_{i=2}^\infty \frac{(\theta U)^i}{i!}\leq 
\frac{\theta^2}{2}\sum_{i=2}^\infty \frac{(\theta U)^{i-2}}{3^{i-2}}\leq \frac{\theta^2}{2}\frac{1}{1-\theta U/3}.
\end{align}
Moreover, 
\begin{align}
\label{eq:n30}
&
\frac{\exp(\theta t)}{\phi(\theta t)}=1+\frac{1+\theta t}{e^{\theta t} - \theta t - 1}\leq 1+\frac{1+\theta t}{\l( \theta t \r)^2/2+\l( \theta t \r)^3/6}\leq 1+\frac{6}{\l( \theta t \r)^2}.
\end{align}
Combination of \eqref{eq:n10},\eqref{eq:n20},\eqref{eq:n30} implies that for $\theta<3/U$, 
\[
\prob{\lmax{S_n}>t}\leq \tr\left(\frac{\mb E S_n^2}{\sigma^2}\right) 
\exp\l[\frac{\theta^2}{2}\frac{1}{1-\theta U/3}\sigma^2 - \theta t\r] \l( 1+\frac{6}{\l( \theta t \r)^2} \r).
\]
Finally, setting $\theta:=\frac{t}{\sigma^2 + Ut/3}$ yields 
\begin{align}
\label{eq:ineq-full}
\prob{\lmax{S_n}>t}\leq \tr\left(\frac{\mb E S_n^2}{\sigma^2}\right) 
\exp\l[ -\frac{t^2/2}{\sigma^2 + Ut/3} \r]
\l( 1+\frac{6}{\l( \theta t \r)^2} \r). 
\end{align}
This inequality is typically interesting only for sufficiently large values of $t$, for instance such that 
$t^2\geq \sigma^2 + Ut/3,$ or $t\geq \frac{1}{6}\l( U+\sqrt{U^2 + 36\sigma^2} \r)$. 
In this range of $t$, $1 + \frac{\l(\sigma^2 + Ut/3 \r)^2}{t^4}\leq 7$ so the inequality simplifies to 
\[
\prob{\left\|\sum_{i=1}^n X_i\right\|>t}\leq 7 \, \tr\left(\frac{\mb E S_n^2}{\sigma^2}\right) 
\exp\left[- \frac{t^2/2}{\sigma^2+tU/3} \right].
\]
It remains to repeat the argument with $X_i$'s replaced by $(-X_i)$'s to obtain an identical bound 
for 
\[
\prob{\lmax{-S_n}>t} = \prob{\lmin{S_n}<-t},
\] 
which yields the desired inequality for the operator norm.
\end{proof}
Next, we state the inequality for general rectangular matrices. 

\begin{corollary}
\label{cor:1}
Let $X_1,\ldots, X_n\in \mb C^{d_1\times d_2}$ be a sequence of independent random matrices such that $\mb EX_i=0, \ i=1,\ldots,n$ and 
$\sigma^2\geq \max\l( \l\|\sum\limits_{i=1}^n \mb E X_iX_i^\ast\r\|, \l\|\sum\limits_{i=1}^n \mb E X_i^\ast X_i\r\| \r)$. 
Moreover, assume that $\|X_i\|\leq U$ almost surely for all $1\leq i\leq n$ and some positive $U\in \mb R$. 
Then, for any $t\geq \frac{1}{6}\l( U+\sqrt{U^2 + 36\sigma^2} \r)$, 
\begin{align}
\label{eq:a05}
&
\prob{\left\|\sum_{i=1}^n X_i\right\|>t}\leq 28 \tilde d
\exp\left[- \frac{t^2/2}{\sigma^2+tU/3} \right],
\end{align}
where $\tilde d = \max\l( r\left(\sum\limits_{i=1}^n \mb EX_i X_i^\ast\right), r\left(\sum\limits_{i=1}^n \mb EX_i^\ast X_i\right)  \r)$.
\end{corollary}
\begin{proof}
We will apply Theorem \ref{th:independent} to the sequence $\m H(X_1),\ldots,\m H(X_n)$ of Hermitian dilations defined in \eqref{eq:dilation}. 
To this end, note that 
\begin{align*}
\mb E\l( \sum_{i=1}^n \m H(X_i)^2 \r)= 
\begin{pmatrix}
\sum_{i=1}^n \mb E X_i X_i^\ast & 0 \\
0 & \sum_{i=1}^n \mb E X_i^\ast X_i
\end{pmatrix},
\end{align*}
hence $\l\| \mb E\l( \sum_{i=1}^n \m H^2(X_i) \r) \r\| = \max\l( \l\|\sum\limits_{i=1}^n \mb E X_iX_i^\ast\r\|, \l\|\sum\limits_{i=1}^n \mb E X_i^\ast X_i\r\| \r)$ and 
the effective rank satisfies 
\[
r\l(\mb E \sum_{i=1}^n \m H^2(X_i) \r) = \frac{\tr(\sum_{i=1}^n \mb E X_i X_i^\ast) + 
\tr(\sum_{i=1}^n \mb E X_i^\ast X_i)}{\max\l( \l\|\sum\limits_{i=1}^n \mb E X_iX_i^\ast\r\|, \l\|\sum\limits_{i=1}^n \mb E X_i^\ast X_i\r\| \r)}
\leq 
2 \max\l( r\left(\sum\limits_{i=1}^n \mb EX_i X_i^\ast\right), r\left(\sum\limits_{i=1}^n \mb EX_i^\ast X_i\right)  \r).
\]
Since $\l\| \sum_{i=1}^n \m H(X_i) \r\| = \l\| \sum_{i=1}^n X_i \r\|$, we get that 
\[
\prob{\left\|\sum_{i=1}^n X_i\right\|>t}\leq 14 \, r\l( \mb E \l[ \sum_{i=1}^n \m H^2(X_i) \r] \r)
\exp\left[- \frac{t^2/2}{\sigma^2+tU/3} \right],
\]
and result follows.
\end{proof}
The bound for $\mb E \l\| \sum_{i=1}^n X_i \r\|$ can be obtained using the tools developed in the course of the proof of Theorem \ref{th:independent}: 
\begin{corollary}
\label{cor:2}
Let the assumptions of corollary \ref{cor:1} be satisfied. Then
\[
\mb E \l\| \sum_{i=1}^n X_i \r\| \leq \frac{5}{4}
\max\l( 2  \sigma \log^{1/2}\l( 2+ 4 \tilde d \r), U\log\l( 2 + 4 \tilde d \r)\r).
\]
where $\tilde d = \max\l( r\left(\sum\limits_{i=1}^n \mb EX_i X_i^\ast\right), r\left(\sum\limits_{i=1}^n \mb EX_i^\ast X_i\right)  \r)$. 
\end{corollary}
\begin{proof}
Let $\phi(x)=e^x - x - 1$ and $\theta>0$. 
Since $\phi$ is a convex function,  
\begin{align}
\label{eq:nn10}
&
\phi\l( \mb E \,\theta \, \l\| \sum_{i=1}^n \m H(X_i) \r\| \r) \leq 
\mb E \phi\l( \theta \l\| \sum_{i=1}^n \m H(X_i) \r\| \r) \leq 
r\l( \mb E \l[ \sum_{i=1}^n \m H^2(X_i) \r] \r) \exp\l( \frac{\theta^2}{2}\frac{\sigma^2}{1-\theta U/3}\r),
\end{align}
where the last inequality follows from \eqref{eq:n1}, \eqref{eq:n2} and \eqref{eq:n20}. 
Next, it is easy to see that $\phi(x)\geq \frac{e^x}{2} - 1$ for all $x>0$, hence \eqref{eq:nn10} implies that
\begin{align*}
e^{\theta \mb E \l\| \sum_{i=1}^n X_i \r\|} & \leq 
2+ 2 r\l(\mb E\sum_{i=1}^n \m H^2(X_i) \r) \exp\l( \frac{\theta^2}{2}\frac{\sigma^2}{1-\theta U/3}\r) \\
&\leq 
\l( 2 + 2 r\l( \mb E\sum_{i=1}^n \m H^2(X_i) \r) \r)\exp\l( \frac{\theta^2}{2}\frac{\sigma^2}{1-\theta U/3}\r),
\end{align*}
or 
$\mb E \l\| \sum_{i=1}^n X_i \r\|\leq \inf_{0<\theta<3/U} \l[ \frac{\log\l( 2(1+ r\l( \mb E\sum_{i=1}^n \m H^2(X_i) \r) \r)}{\theta}
+\frac{\theta}{2} \frac{\sigma^2}{1-\theta U/3} \r]$. 
Consider two cases: 
\begin{enumerate}[(a)]
\item 
If $\frac{\log^{1/2}\l( 2(1+ r\l( \mb E\sum_{i=1}^n \m H^2(X_i) \r) \r)}{\sigma}\leq \frac{2}{U}$, setting 
$\theta:=\frac{\log^{1/2}\l( 2(1+ r\l( \mb E\sum_{i=1}^n \m H^2(X_i) \r) \r)}{\sigma}$ gives 
\[
\mb E \l\| \sum_{i=1}^n X_i \r\|\leq \frac{5}{2} \sigma \log^{1/2}\l( 2(1+ r\l( \mb E\sum_{i=1}^n \m H^2(X_i) \r) \r).
\]
\item 
If $\frac{\log^{1/2}\l( 2(1+ r\l( \mb E\sum_{i=1}^n \m H^2(X_i) \r) \r)}{\sigma} > \frac{2}{U}$, 
set $\theta:=\frac{2}{U}$ to get 
\[
\mb E \l\| \sum_{i=1}^n X_i \r\|\leq \frac{5}{4} U \log\l( 2(1+ r\l( \mb E\sum_{i=1}^n \m H^2(X_i) \r) \r),
\]
\end{enumerate}
and the claim follows since $ r\Big( \mb E\sum_{i=1}^n \m H^2(X_i) \Big) \leq 2\tilde d$.  
\end{proof}


\subsection{Concentration inequality for the sums of martingale differences}

Our next goal is to obtain a version of Friedman's inequality \cite{freedman1975tail} for the sums of matrix-valued martingale differences (see chapter 8 of \cite{pollard2012convergence} for the definitions). 
For $t\in \mb R$, define $p(t):=\min(-t,1)$. 
Note that 
\begin{enumerate}[(a)]
\item $p(t)$ is concave;
\item $g(t):=e^t-1+p(t)$ is non-negative for all $t$ and increasing for $t>0$.
\end{enumerate}
In particular, since $p(t)$ is concave, it follows from fact \ref{fact:06} and Jensen's inequality that for any random self-adjoint matrix $Y$ such that $\mb EY$ is well-defined,
\begin{align}
\label{z4}
&
\mb E\, \tr \, p(Y)\leq \tr\, p(\mb EY).
\end{align}
We are ready to state the main result of this section.
\begin{theorem}
\label{th:martingale}
Let $X_1,\ldots, X_n$ be a sequence of martingale differences with values in the set of $d\times d$ self-adjoint matrices and such that $\|X_i\|\leq U$ almost surely for all $1\leq i\leq n$ and some positive $U\in \mb R$. 
Denote $W_n:=\sum\limits_{i=1}^n \mb E_{i-1} X_i^2$. 
Then for any $t\geq \frac{1}{6}\l( U+\sqrt{U^2 + 36\sigma^2} \r)$,
\begin{align*}
&
\prob{\left\|\sum_{i=1}^n X_i\right\|>t, \ \lmax{W_n}\leq \sigma^2}\leq 
50 \,\tr \l[ p\l( - \frac{t}{U} \frac{\mb EW_n}{\sigma^2} \r) \r] \exp\l( -\frac{t^2/2}{\sigma^2 + tU/3} \r).
\end{align*}
\end{theorem}

\begin{remark}
Expression $\tr \l[ p\l( - \frac{t}{U} \frac{\mb EW_n}{\sigma^2} \r) \r]$ which plays the role of the ``dimension factor'' in our bound has a very simple meaning: acting on the cone of nonnegative definite operators, the function $A\mapsto p(-A)$ just truncates the eigenvalues of $A$ on the unit level. 
It is easy to see that if the eigenvalues of $\mb EW_n$ decay polynomially, that is, $\lambda_j \l( \mb EW_n \r)\leq C_1 \frac{\sigma^2}{j^r}$ for some $C_1>0$ and $r>1$, then 
\[
\tr \l[ p\l( - \frac{t}{U} \frac{\mb EW_n}{\sigma^2} \r) \r] \leq \min(d, c_2(t/U)^{1/r}).
\]
In particular, this gives an improvement over the multiplicative factor $t$ appearing in the bounds established in \cite{Hsu2012Tail-inequaliti00}.
\end{remark}
\begin{proof}
Recall that $\phi(\theta)=e^{\theta}-\theta-1$, and set $S_n:=\sum\limits_{i=1}^n X_i$. 
Let $t>0$, and assume $\theta$ is such that $\theta t-\frac{\phi(\theta U)}{U^2}\sigma^2>0$. 
Define event $\m E$ via
\[
\m E:=\left\{\lmax{\theta S_n - \frac{\phi(\theta U)}{U^2} W_n}\geq \theta t - \frac{\phi(\theta U)}{U^2}\sigma^2\right\}.
\]
Note that by the triangle inequality 
$
\m E\supseteq\left\{\lmax{S_n}\geq t, \ \lmax{W_n}\leq \sigma^2\right\}.
$ 
We proceed by estimating $\prob{\m E}$ in a way similar to \eqref{z2}:
\begin{align}
\label{z5}
\nonumber
\prob{\m E}&
= \prob{\lmax{g\l( \theta S_n - \frac{\phi(\theta U)}{U^2} W_n \r)}\geq g \l( \theta t - \frac{\phi(\theta U)}{U^2}\sigma^2\r)} \\
&
\leq
\frac{ \tr\,\mb E\l(g\l( \theta S_n - \frac{\phi(\theta U)}{U^2} W_n \r)\r) }{g\l( \theta t - \frac{\phi(\theta U)}{U^2}\sigma^2 \r)}.
\end{align}
Let 
$
Y_k:=\tr\exp(\theta S_k - \frac{\phi(\theta U)}{U^2} W_k). 
$ 
It follows from lemma 2.1 in  \cite{tropp2} that $Y_k$ is a supermartingale with initial value $d$, hence 
$
\mb E\,\tr \, \exp\l(\theta S_n - \frac{\phi(\theta U)}{U^2}W_n\r)\leq d.
$
Together with (\ref{z4}), this gives
\begin{align}
\label{n4}
\nonumber
\tr\,\mb E\,g\l( \theta S_n - \frac{\phi(\theta U)}{U^2} W_n \r)&
=\tr \, \mb E\l(\exp\l( \theta S_n - \frac{\phi(\theta U)}{U^2}W_n \r)-I_d+p\l( \theta S_n - \frac{\phi(\theta U)}{U^2}W_n \r)\r) \\
\nonumber
&
\leq \mb E\,\tr\, p\l(\theta S_n - \frac{\phi(\theta U)}{U^2} W_n\r) 
\leq 
\tr\,p\l(\theta\mb E S_n - \frac{\phi(\theta U)}{U^2}W_n\r) \\ 
&
= \tr\,p\l(- \frac{\phi(\theta U)}{U^2}\mb E W_n\r).
\end{align}
Set $\theta:=\frac{1}{U}\log\l(1 + Ut/\sigma^2 \r)$. 
Since $\mb EW_n \succeq 0$ and
$
\phi(y)\leq e^{y}-1
$ 
for $y\geq 0$, $\frac{\phi(\theta U)}{U^2}\mb E W_n \preceq \frac{e^{\theta U}-1}{U^2}\mb EW_n$, hence 
\[
\tr p\l(- \frac{\phi(\theta U)}{U^2}\mb E W_n\r)\leq \tr p\l(- \frac{e^{\theta U}-1}{U^2}\mb EW_n \r) = 
\tr p\l( - \frac{t}{U} \frac{\mb EW_n}{\sigma^2} \r)
\]
by monotonicity of $p(-t)=\min(t,1)$.  
Finally, for our choice of $\theta$,
\begin{align}
\nonumber
\l( g\l( \theta t - \frac{\phi(\theta U)}{U^2}\sigma^2 \r) \r)^{-1} & = 
e^{-\theta t + \frac{\phi(\theta U}{U^2}\sigma^2} 
\frac{e^{\theta t - \frac{\phi(\theta U)}{U^2}\sigma^2}}{g\l( \theta t - \frac{\phi(\theta U))}{U^2}\sigma^2 \r)} 
\\
&
\label{eq:nn20}
= \exp\l(-\frac{\sigma^2}{U^2} h\l(Ut/\sigma^2 \r) \r) 
\frac{ \exp\l(\frac{\sigma^2}{U^2} h\l(Ut/\sigma^2 \r) \r)  }
{ \exp\l(\frac{\sigma^2}{U^2} h\l(Ut/\sigma^2 \r) \r) - \frac{\sigma^2}{U^2} h\l(Ut/\sigma^2 \r) - 1 },
\end{align}
where $h(z):=(1+z)\log(1+z) - z$ for $z\geq 0$. 
A well-known numerical inequality states that 
\begin{align}
\label{eq:nn30}
h(z)\geq \frac{z^2}{2(1+z/3)},
\end{align} 
hence the first term of the product in \eqref{eq:nn20} satisfies 
\[
\exp\l(-\frac{\sigma^2}{U^2} h\l(Ut/\sigma^2 \r) \r) \leq \exp\l( -\frac{t^2/2}{\sigma^2 + tU/3} \r). 
\] 
By \eqref{eq:n30} and \eqref{eq:nn30}, the second term in \eqref{eq:nn20} can be estimated as  
\begin{align*}
\frac{ \exp\l(\frac{\sigma^2}{U^2} h\l(Ut/\sigma^2 \r) \r)  }
{ \exp\l(\frac{\sigma^2}{U^2} h\l(Ut/\sigma^2 \r) \r) - \frac{\sigma^2}{U^2} h\l(Ut/\sigma^2 \r) - 1 }
\leq 1 + \frac{6}{\l(\frac{\sigma^2}{U^2}h(Ut/\sigma^2) \r)^2}
\leq 
1+ \frac{6}{\l(\frac{t^2/2}{\sigma^2 +Ut/3} \r)^2}.
\end{align*}
Whenever $t^2\geq \sigma^2 +Ut/3$, the latter expression is further bounded by $25$, which implies that in this range of values of $t$, 
\[
\prob{\m E} \leq 25 \,\tr p\l( - \frac{t}{U} \frac{\mb EW_n}{\sigma^2} \r) \exp\l( -\frac{t^2/2}{\sigma^2 + tU/3} \r).
\]
Repeating the argument with $X_i$'s replaced by $(-X_i)$'s concludes the proof.
\end{proof}
We note that Theorem \ref{th:martingale} admits generalization to the case of rectangular matrices, in a way that is very similar to Corollary \ref{cor:1}. We omit the details.

\subsection{Extensions to the case of self-adjoint operators}
\label{sec:extension}

Unlike the results which depend on the dimension of the ambient space, Theorems \ref{th:independent} and \ref{th:martingale} apply to the case when $X_1,\ldots,X_n$ is a sequence of self-adjoint Hilbert-Schmidt operators $X_i: \mb H\mapsto \mb H$ acting on a separable Hilbert space $\l(\mb H, \dotp{\cdot}{\cdot}_{\mb H}\r)$, 
such that $\ker(\mb EX_i)=\mb H$, $i=1,\ldots,n$. 
Here, $\mb EX$ is an operator such that $\dotp{\mb (\mb EX) z_1}{z_2}_{\mb H}=\mb E \dotp{X z_1}{z_2}_\mb H$ for any $z_1,z_2\in \mb H$.
We will formally show how to extend results of Theorem \ref{th:independent}; similar argument applies to Theorem \ref{th:martingale}.

The following argument relies on the fact that a compact operator can be approximated by a sequence of operators with finite rank. 
Let $L_1\subset L_2\subset\ldots$ be a nested sequence of finite dimensional subspaces of $\mb H$ such that 
$\bigcup\limits_j L_j=\mb H$, and let 
$P_{L_j}$ and $P_{L_j^\perp}$ be the orthogonal projectors onto $L_j$ and its orthogonal complement $L_j^\perp$ respectively. 
For any fixed $j$, we will apply Theorem \ref{th:independent} to a sequence of finite dimensional operators $\l\{P_{L_j}X_i P_{L_j}\r\}_{i\geq 1}$ mapping $L_j$ to itself. 
Note that for any $v\in \mb H$, $j\geq 1$ and $1\leq i\leq n$, 
\[
\l\| X_i v - P_{L_j}X_i P_{L_j} v \r\|_\mb H \leq \l\| X_i P_{L_j^\perp} v \r\|_\mb H  + \l\| P_{L_j^\perp} X P_{L_j} v  \r\|_\mb H,
\] 
hence 
$\l\|\sum\limits_{i=1}^n \l(X_i - P_{L_j}X_i P_{L_j}\r)\r\|\xrightarrow{} 0$ almost surely as $j\to\infty$, and
\begin{align}
\label{eq:a10}
&
\prob{\left\|\sum_{i=1}^n X_i\right\|>t}
\leq \liminf_{j\to\infty}\prob{\l\|\sum\limits_{i=1}^n P_{L_j}X_i P_{L_j} \r\|>t}
\end{align}
by Fatou's lemma. 
Since $A\preceq B$ implies $SAS^*\preceq SBS^*$ (fact \ref{fact:04}), taking $A=P^2_{L_j}=P_{L_j}, \ B=I$ and $S=P_{L_j}X$ gives
$
(P_{L_j}XP_{L_j})^2\preceq P_{L_j}X^2 P_{L_j}
$,
thus 
\begin{align}
\label{eq:a15}
\liminf_{j\to\infty}&
\frac{\tr\left(\sum\limits_{i=1}^n \mb E\l(P_{L_j}X_i P_{L_j}\r)^2\r)}{\l\| \sum\limits_{i=1}^n \mb E (P_{L_j}X_i P_{L_j})^2 \r\| }
\leq
\liminf_{j\to\infty}\frac{\tr\l(\sum\limits_{i=1}^n \mb E\l(P_{L_j}X_i^2P_{L_j}\r)\r)}{\l\| \sum\limits_{i=1}^n \mb E (P_{L_j}X_iP_{L_j})^2 \r\|}
 \\
&\nonumber
\leq 
\frac{\tr\left(\sum\limits_{i=1}^n \mb EX_i^2\right)}{\limsup\limits_{j\to\infty} \l\| \sum\limits_{i=1}^n \mb E(P_{L_j}X_i P_{L_j})^2 \r\|}=
 \frac{\tr\left(\sum\limits_{i=1}^n \mb EX_i^2\right)}{\l\| \sum\limits_{i=1}^n \mb EX_i^2\r\|},
\end{align}
where in the last step we used a simple bound 
\begin{align}
\label{eq:a20}
\l\|X_i^2-(P_{L_j}X_iP_{L_j})^2\r\|&
\nonumber
=\l\|X_i (X_i - P_{L_j}X_iP_{L_j})+(X_i-P_{L_j}X_iP_{L_j})P_{L_j}X_iP_{L_j}\r\|\leq \\
&
\leq 2\|X_i\| \l\| X_i-P_{L_j}X_iP_{L_j} \r\|
\end{align}
which converges to $0$ almost surely. 
Since $\|X_i\|\leq U$ a.s., (\ref{eq:a20}) implies by the dominated convergence theorem that 
\begin{align}
\label{eq:a25}
&
\l\| \sum\limits_{i=1}^n \mb E \l( P_{L_j}X_i P_{L_j} \r)^2 - \sum\limits_{i=1}^n\mb E X_i^2 \r\| \leq 
\mb E\l\|\sum\limits_{i=1}^n(P_{L_j}X_i P_{L_j})^2 - \sum\limits_{i=1}^n X_i^2 \r\|\to 0 \text{ as } j\to\infty.
\end{align}
It remains to apply Theorem \ref{th:independent} to the right-hand side of (\ref{eq:a10}). 
Combined with (\ref{eq:a15}) and (\ref{eq:a25}), it implies that 
\[
\prob{\left\|\sum_{i=1}^n X_i\right\|>t}\leq 14 \, r\left(\sum\limits_{i=1}^n \mb EX_i^2\right)
\exp\left[- \frac{t^2/2}{\sigma^2+tU/3} \right]
\]
whenever $t\geq \frac{1}{6}\l( U+\sqrt{U^2 + 36\sigma^2} \r)$.

\section{Example: vector-valued Bernstein's inequality}
\label{sec:applications}
Several relevant applications of the bounds that depend on the effective rank, such as covariance estimation and approximate matrix multiplication, were demonstrated in \cite{Hsu2012Tail-inequaliti00}. 
Our bounds apply in those examples as well and yield slightly sharper tails. 
We demonstrate another immediate corollary of our results. 

\begin{corollary}
\label{corr:bernstein}
Let $Y_1,\ldots,Y_n\in \mb C^d$ be a sequence of independent random vectors such that $\mb EY_i=0$ and $\|Y_i\|_2\leq U$ almost surely for all $1\leq i\leq n$ and some $U>0$. 
Denote $\sigma^2:=\sum\limits_{i=1}^n \mb E\l\| Y_i \r\|_2^2$.
Then for all $t\geq \frac{1}{6}\l( U+\sqrt{U^2 + 36\sigma^2} \r)$
\[
\prob{\l\|\sum\limits_{i=1}^n Y_i\r\|_2>t}\leq 28 \exp\left[- \frac{t^2/2}{\sigma^2+tU/3} \right].
\]
\end{corollary}
\begin{remark}
This result can be compared to the bound in \cite{Ledoux1991Probability-in-00} (see formula 6.13) obtained by a combination of classical martingale methods and a trick of V. Yurinskii; also, see \cite{Ahmad2013Probability-ine00} for a different version of Bernstein's inequality.
\end{remark}
\begin{proof}
We proceed as in the proof of corollary \ref{cor:1}. 
Observe that 
$\m H(Y_i)^2=\begin{pmatrix}
\|Y_i\|_2^2 & 0 \\
0 & Y_i Y_i^*
\end{pmatrix}
$, where $\m H(\cdot)$ is the self-adjoint dilation defined in \eqref{eq:dilation}. 
Clearly, 
\[
\l\|\sum\limits_{i=1}^n \mb E \m H^2(Y_i)\r\|=
\l\|
\begin{pmatrix}
\sum\limits_{i=1}^n\mb E\|Y_i\|_2^2 & 0 \\
0 & \sum\limits_{i=1}^n \mb EY_i Y_i^\ast
\end{pmatrix}
\r\|
=\sum\limits_{i=1}^n\mb E\|Y_i\|_2^2
\]
and, similarly, $\tr\l(\sum\limits_{i=1}^n \mb E \m H^2(Y_i)\r) = 2\sum\limits_{i=1}^n \mb E \|Y_i\|_2^2$, hence the effective rank satisfies
$r\l( \sum\limits_{i=1}^n \mb E \m H^2(Y_i) \r)=2$. 
The result now follows by repeating the proof of corollary \ref{cor:1}.
\end{proof}

\bibliographystyle{alpha}
\bibliography{bibliography}

\newcommand{\etalchar}[1]{$^{#1}$}
\begin{thebibliography}{KL{\etalchar{+}}16}

\bibitem[AA13]{Ahmad2013Probability-ine00}
I.~A. Ahmad and M.~Amezziane.
\newblock Probability inequalities for bounded random vectors.
\newblock {\em Statist. Probab. Lett.}, 83(4):1136--1142, 2013.

\bibitem[AW02]{Ahlswede2002Strong-converse00}
R.~Ahlswede and A.~Winter.
\newblock Strong converse for identification via quantum channels.
\newblock {\em IEEE Trans. Inform. Theory}, 48(3):569--579, 2002.

\bibitem[Bha97]{bhatia1997matrix}
R.~Bhatia.
\newblock Matrix analysis.
\newblock {\em Springer}, 1997.

\bibitem[Car10]{carlen2010trace}
E.~Carlen.
\newblock Trace inequalities and quantum entropy: an introductory course.
\newblock {A}vailable at
  \url{http://www.mathphys.org/AZschool/material/AZ09-carlen.pdf}, 2010.

\bibitem[Fre75]{freedman1975tail}
David~A Freedman.
\newblock On tail probabilities for martingales.
\newblock {\em the Annals of Probability}, pages 100--118, 1975.

\bibitem[Gol65]{sidne1965lower}
S.~Golden.
\newblock Lower bounds for the {H}elmholtz function.
\newblock {\em Phys. Rev. (2)}, 137:B1127--B1128, 1965.

\bibitem[HKZ12]{Hsu2012Tail-inequaliti00}
D.~Hsu, S.~M. Kakade, and T.~Zhang.
\newblock Tail inequalities for sums of random matrices that depend on the
  intrinsic dimension.
\newblock {\em Electronic Communications in Probability}, 17(14):1--13, 2012.

\bibitem[KL{\etalchar{+}}16]{koltchinskii2016asymptotics}
Vladimir Koltchinskii, Karim Lounici, et~al.
\newblock Asymptotics and concentration bounds for bilinear forms of spectral
  projectors of sample covariance.
\newblock In {\em Annales de l'Institut Henri Poincar{\'e}, Probabilit{\'e}s et
  Statistiques}, volume~52, pages 1976--2013. Institut Henri Poincar{\'e},
  2016.

\bibitem[Lie73]{lieb1}
E.~H. Lieb.
\newblock Convex trace functions and the {W}igner-{Y}anase-{D}yson conjecture.
\newblock {\em Advances in Math.}, 11:267--288, 1973.

\bibitem[LP86]{lust1986inegalites}
Fran{\c{c}}oise Lust-Piquard.
\newblock In\'{e}galit\'{e}s de {Khintchine} dans ${C}_p$ ($1<p<\infty$).
\newblock {\em CR Acad. Sci. Paris}, 303:289--292, 1986.

\bibitem[LPP91]{lust1991non}
Fran{\c{c}}oise Lust-Piquard and Gilles Pisier.
\newblock Non commutative khintchine and paley inequalities.
\newblock {\em Arkiv f{\"o}r matematik}, 29(1):241--260, 1991.

\bibitem[LT91]{Ledoux1991Probability-in-00}
M.~Ledoux and M.~Talagrand.
\newblock {\em Probability in {B}anach spaces}, volume~23.
\newblock Springer-Verlag, Berlin, 1991.
\newblock Isoperimetry and processes.

\bibitem[Oli10]{Oliveira2010Concentration-o00}
R.~I. Oliveira.
\newblock Concentration of the adjacency matrix and of the {L}aplacian in
  random graphs with independent edges.
\newblock {\em arXiv preprint arXiv:0911.0600}, 2010.

\bibitem[Pol12]{pollard2012convergence}
David Pollard.
\newblock {\em Convergence of stochastic processes}.
\newblock Springer Science \& Business Media, 2012.

\bibitem[Rud99]{rudelson1999random}
Mark Rudelson.
\newblock Random vectors in the isotropic position.
\newblock {\em Journal of Functional Analysis}, 164(1):60--72, 1999.

\bibitem[Tho65]{thompson1965inequality}
C.~J. Thompson.
\newblock Inequality with applications in statistical mechanics.
\newblock {\em J. Mathematical Phys.}, 6:1812--1813, 1965.

\bibitem[Tro11]{tropp2}
J.~A. Tropp.
\newblock Freedman's inequality for matrix martingales.
\newblock {\em Electron. Commun. Probab.}, 16:262--270, 2011.

\bibitem[Tro12]{tropp2012user}
J.~A. Tropp.
\newblock User-friendly tail bounds for sums of random matrices.
\newblock {\em Foundations of computational mathematics}, 12(4):389--434, 2012.

\end{thebibliography}

\appendix
\section{Proof of Lemma \ref{lemma:semidef}}
\label{sec:lemma_pf}
\noindent Writing the series expansion for $e^{\theta X_i}$ and using the fact that $\mb EX_i=0$, we obtain
\begin{align*}
\mb Ee^{\theta X_i}& = I + \mb E\l[ \frac{\theta^2 X_i^2}{2!}+\frac{\theta^3 X_i^3}{3!}+\ldots+\frac{(\theta X_i)^{k}}{k!}+\ldots \r]  \\
&
= I + \theta^2 \mb E\l[ X_i\l(\frac{1}{2!}+\frac{\theta X_i}{3!}+\ldots+\frac{(\theta X_i)^{k-1}}{(k+1)!}+\ldots\r) X_i\r] \\
&
\preceq 
I + \theta^2\mb E \l[ X_i \l(\frac{1}{2!}+\frac{\theta \| X_i \|}{3!}+\ldots+\frac{\theta^k \|X_i\|^{k}}{(k+1)!}+\ldots\r) X_i \r] \\
&
= I +\theta^2 \mb E \l[ X_i^2 \l(\frac{e^{\theta \|X_i\|}-\theta\|X_i\|-1}{\theta^2\|X_i\|^2}\r) \r] 
\preceq I + \frac{ \phi(\theta U)}{U^2} \mb E X_i^2.
\end{align*}
We used the semidefinite relation 
$\l(\frac{1}{2!}+\ldots+\frac{(\theta X_i)^{k-1}}{(k+1)!}+\ldots\r) \preceq  I \l(\frac{1}{2!}+\frac{\theta \| X_i \|}{3!}+\ldots+\frac{\theta^{k-1} \|X_i\|^{k}}{(k+1)!}+\ldots\r)$ 
and fact \ref{fact:04} to move from the second to the third line,
and the assumption $\|X_i\|\leq U$ together with monotonicity of $g(s) = \frac{e^s - s - 1}{s^2}$ to obtain the last semidefinite relation. 

Now \eqref{eq:z1} follows from the relation $I + \frac{\phi(\theta U)}{U^2} \mb E X_i^2 \preceq e^{\frac{\phi(\theta U)}{U^2} \mb E X_i^2}$ which is a consequence of fact \ref{fact:01}, the elementary inequality $1+x\leq e^x$, and the operator monotonicity of the logarithm (fact \ref{fact:03}).






\end{document}